\documentclass{scrartcl}

\usepackage[english]{babel}
\usepackage{url}
\usepackage{latexsym}
\usepackage{amsmath}
\usepackage{amsthm}
\usepackage{stmaryrd}
\usepackage[linesnumbered,ruled]{algorithm2e}

\begin{document}
\title{Examples of distance magic labelings of the $6$-dimensional hypercube}
\author{Petr Savick\'y\footnote{Institute of Computer Science of the Czech
Academy of Sciences, Czech Republic, savicky@cs.cas.cz}}
\date{\empty}
\maketitle

\begin{abstract}
A distance magic labeling of an $n$-dimensional hypercube is
a labeling of its vertices by natural numbers from $\{0, \ldots, 2^n-1\}$, such
that for all vertices $v$ the sum of the labels of the neighbors of $v$
is the same. Such a labeling is called neighbor-balanced, if, moreover, for each
vertex $v$ and an index $i=0,\ldots,n-1$, exactly half of the neighbors of
$v$ have digit $1$ at $i$-th position of the binary representation of their
label. We demonstrate examples of non-neighbor-balanced distance magic labelings
of $6$-dimensional hypercube obtained by a SAT solver.
\end{abstract}

\section{Introduction}
By the results of \cite{GK13}, a distance magic labeling (DML) of
an $n$-dimensional hypercube exists if and only if $n \equiv 2 \ (\mathrm{mod}\ 4)$.
One can easily verify that if such a labeling exists, then for each vertex,
the sum of the labels of its neighbors is $\frac{1}{2}n(2^n-1)$.
In order to prove the existence of the required labelings, the authors provide
a construction of a neighbor-balanced
distance magic labeling for each $n$ satisfying the condition above
disproving a conjecture formulated in \cite{ARSP04}.
In this note, we present a non-neighbor-balanced distance magic labeling
of the $6$-dimensional hypercube $Q_6$ which solves Problem 3.7 formulated
in \cite{GK13}.

\section{SAT instance}

In order to encode the search for distance magic labelings of $Q_n$
into a SAT instance, we use $n2^n$ variables to encode
binary digits of the labels. In order to enforce the required conditions
on the labels we use the following simple approach with $n=6$.
The condition that two $n$-digit binary numbers are different is
encoded using a formula with $n$ auxiliary variables and $2n+1$ clauses.
The condition that the labels are pairwise different is expressed
using ${2^n \choose 2}$ copies of this formula with disjoint sets of
auxiliary variables. The condition that the
sum of $n$ binary $n$-digit numbers is $\frac{1}{2}n(2^n-1)$ can be
expressed in different ways. For $n=6$, $\frac{1}{2}n(2^n-1)=189$
and the formula was obtained using PySAT \cite{PySAT} interface
to PBLib \cite{PBLib}.
A copy of this formula with new auxiliary variables is used to express
the sum condition on the labels of each vertex in $Q_6$.
Additionally, the formula was extended with literals enforcing
the label of the zero vertex and the labels of its neighbors. The
resulting formula
consists of $26176$ variables and $67146$ clauses and can be solved
by CaDiCal SAT solver \cite{CaDiCal} in a few minutes. For $20$
different seeds, the running time on Intel Xeon processor
(X5680, 3.33 GHz) was between $13$ sec and $768$ sec with
average $175$ sec. For each of these seeds, the solver
obtained a different solution and the first $5$ of them
are presented below.

\section{Examples of non-neighbor-balanced DML of $Q_6$}

The formula used to construct the labelings enforces that the
zero vertex has label $0$ and its neighbors have labels
$4, 6, 36, 38, 52, 53$. The binary expansions of these numbers are
$$
\begin{array}{r|rrrrrr}
 4 & 0 & 0 & 0 & 1 & 0 & 0\\
 6 & 0 & 0 & 0 & 1 & 1 & 0\\
36 & 1 & 0 & 0 & 1 & 0 & 0\\
38 & 1 & 0 & 0 & 1 & 1 & 0\\
52 & 1 & 1 & 0 & 1 & 0 & 0\\
53 & 1 & 1 & 0 & 1 & 0 & 1\\
\end{array}
$$
so all the obtained labelings are non-neighbor-balanced. Verification
of correctness of each of the labelings can be done by computing
the $64$ required sums and is left to the reader.

Examples of the obtained labelings of $Q_6$ are presented in Table \ref{tab1}
in the appendix as tables $8 \times 8$
assuming that $\{0,1\}^6 = \{0,1\}^3 \times \{0,1\}^3$ where the rows
of the tables
correspond to all possible combinations of the first three components
in the lexicographic order and, similarly, the columns correspond
to the combinations of the last three components.

\appendix

\begin{table}[h]
\begin{center}
\begin{tabular}{rrrrrrrr}
 0 &  4 &  6 & 15 & 36 & 60 & 47 & 49 \\
38 &  8 & 61 & 46 & 41 & 45 & 28 & 20 \\
52 & 56 & 29 & 30 & 40 & 24 & 51 &  5 \\
21 & 44 & 32 & 62 &  9 & 26 & 13 & 10 \\
53 & 50 & 37 & 54 &  1 & 31 & 19 & 42 \\
58 & 12 & 39 & 23 & 33 & 34 &  7 & 11 \\
43 & 35 & 18 & 22 & 17 &  2 & 55 & 25 \\
14 & 16 &  3 & 27 & 48 & 57 & 59 & 63
\end{tabular}

\vspace{5mm}
\begin{tabular}{rrrrrrrr}
 0 &  4 &  6 & 21 & 36 & 62 & 58 & 54 \\
38 & 14 & 50 & 55 & 28 & 47 & 12 & 19 \\
52 & 32 & 34 & 24 & 23 & 15 & 61 & 22 \\
56 & 43 & 30 & 45 & 17 & 37 &  3 & 10 \\
53 & 60 & 26 & 46 & 18 & 33 & 20 &  7 \\
41 &  2 & 48 & 40 & 39 & 29 & 31 & 11 \\
44 & 51 & 16 & 35 &  8 & 13 & 49 & 25 \\
 9 &  5 &  1 & 27 & 42 & 57 & 59 & 63 \\
\end{tabular}

\vspace{5mm}
\begin{tabular}{rrrrrrrr}
 0 &  4 &  6 & 26 & 36 & 51 & 34 & 45 \\
38 & 23 &  8 & 30 & 56 & 47 & 46 & 54 \\
52 & 41 & 60 & 42 & 35 & 43 & 24 &  5 \\
44 & 31 & 49 & 50 &  1 &  2 & 15 & 10 \\
53 & 48 & 61 & 62 & 13 & 14 & 32 & 19 \\
58 & 39 & 20 & 28 & 21 &  3 & 22 & 11 \\
 9 & 17 & 16 &  7 & 33 & 55 & 40 & 25 \\
18 & 29 & 12 & 27 & 37 & 57 & 59 & 63 \\
\end{tabular}

\vspace{5mm}
\begin{tabular}{rrrrrrrr}
 0 &  4 &  6 & 41 & 36 & 43 & 61 & 49 \\
38 & 35 &  5 & 29 & 60 & 48 & 42 &  7 \\
52 & 37 & 50 & 46 &  1 & 12 & 47 & 19 \\
45 & 54 & 23 & 39 &  8 & 31 & 30 & 10 \\
53 & 33 & 32 & 55 & 24 & 40 &  9 & 18 \\
44 & 16 & 51 & 62 & 17 & 13 & 26 & 11 \\
56 & 21 & 15 &  3 & 34 & 58 & 28 & 25 \\
14 &  2 & 20 & 27 & 22 & 57 & 59 & 63 \\
\end{tabular}

\vspace{5mm}
\begin{tabular}{rrrrrrrr}
 0 &  4 &  6 & 15 & 36 & 54 & 19 & 23 \\
38 & 61 & 60 & 51 & 26 & 45 & 49 & 17 \\
52 & 16 & 62 & 55 & 58 & 42 & 34 & 28 \\
 7 & 22 & 24 & 31 & 13 & 30 & 20 & 10 \\
53 & 43 & 33 & 50 & 32 & 39 & 41 & 56 \\
35 & 29 & 21 &  5 &  8 &  1 & 47 & 11 \\
46 & 14 & 18 & 37 & 12 &  3 &  2 & 25 \\
40 & 44 &  9 & 27 & 48 & 57 & 59 & 63 \\
\end{tabular}
\end{center}
\caption{Five examples of non-neighbor-balanced DMLs of $Q_6$.} \label{tab1}
\end{table}


\begin{thebibliography}{0}
\bibitem{ARSP04}
B.~D.~Acharya, S.~B.~Rao, T.~Singh, and V.~Parameswaran,
Neighborhood magic graphs,
In National Conference on Graph Theory, Combinatorics and Algorithm, 2004.
\bibitem{GK13}
Petr Gregor,
Petr Kov\'a{\v r},
Distance magic labelings of hypercubes,
Electronic Notes in Discrete Mathematics 40 (2013), 145--149.
\bibitem{PySAT}
PySAT: SAT technology in Python,
\url{https://pysathq.github.io/}
\bibitem{PBLib}
PBLib -- A C++ Toolkit for Encoding Pseudo-Boolean Constraints into CNF,\\
\url{http://tools.computational-logic.org/content/pblib.php}
\bibitem{CaDiCal}
CaDiCaL SAT solver,
\url{http://fmv.jku.at/cadical/}
\end{thebibliography}
\end{document}